\documentclass[12pt]{article}
\usepackage{mathrsfs}
\usepackage{amsfonts}
\usepackage{latexsym}
\date{}

\begin{document}
\title{  Finitistic Dimension of Faithfully Flat Weak Hopf-Galois Extension$^\star$}
\author{{\small Aiping Zhang}\\
{\small School of Mathematics and Statistics,\ Shandong University,\ Weihai, Weihai 264209, \ China }\\
 {\small E-mail: aipingzhang@sdu.edu.cn}
\\
}

\pagenumbering{arabic}

\maketitle

\footnote{ $^\star$Supported by the NSF of China (Grant No.
11601274).}

\begin{center}
 \begin{minipage}{120mm}
   \small\rm
   {\bf  Abstract.}\ \  Let $H$ be a finite-dimensional weak Hopf algebra over a field $k$ and $A/B$ be a right faithfully flat weak $H$-Galois extension. We prove that if the finitistic dimension of $B$ is finite, then it is less than or equal to that of $A$. Moreover, suppose that $H$ is semisimple. If the finitistic dimension conjecture holds,  then the finitistic dimension of $B$ is equal to that of $A$.
\end{minipage}
\end{center}

\begin{center}
  \begin{minipage}{120mm}
\end{minipage}
\end{center}

\begin{center}
  \begin{minipage}{120mm}
\end{minipage}
\end{center}
\vskip -0.3cm

\section {Introduction}

\vskip 0.2in

\ \ \ \ \ \ Throughout the paper, we work over a field $k$, vector spaces, algebras, coalgebras and unadorned $\otimes$ are over $k$. Given an algebra $A$, we denote by
$A$-Mod and $A$-mod the categories of left $A$-modules and of finitely generated left $A$-modules, respectively. For a left $A$-module $M$, we denote by
pd$M$ the projective dimension of $M$.

Weak Hopf algebras are generalizations of ordinary Hopf algebras. Examples of weak Hopf algebras are groupoid algebras, face algebras [7], quantum groupoids [9], generalized Kac algebras [13], and so on. The importance of weak Hopf algebras for quantum field theory and operator algebras is well understood and there is much literature on this subject, see [4,8,15].

Recall from B$\ddot{\rm o}$hm et al.[5] that a weak $k$-bialgebra $H$ is both a $k$-algebra ($m, \mu$) and a $k$-coalgebra $(\Delta, \varepsilon$) such that $\Delta(hk) = \Delta(h)\Delta(k)$, and
$$\Delta^2(1) = 1_1 \otimes 1_21^{'}_1 \otimes 1^{'}_2 = 1_1 \otimes 1^{'}_11_2\otimes 1^{'}_2,$$
$$\varepsilon(hkl) = \varepsilon(hk_1)\varepsilon(k_2l) = \varepsilon(hk_2)\varepsilon(k_1l),$$ for all $h, k, l \in H$, where $1^{'}$ stands for another copy of 1. Here we use
Sweedler$^{'}$s notation for the comultiplication. Namely, for $h\in H, \Delta(h) = h_1\otimes h_2$, where we omit the summation symbol and indices. The maps $\varepsilon_t, \varepsilon_s :H \rightarrow H$ defined by

$$ \varepsilon_t(h) = \varepsilon(1_1h)1_2;\ \varepsilon_s(h) = 1_1\varepsilon(h1_2)$$
are called the target map and source map, and their images $H_t$ and $H_s$ are called the target and source space.

A weak Hopf algebra $H$ is a weak bialgebra together with a $k$-linear map $S: H\rightarrow H$(called the antipode) satisfying
$$S*id_H = \varepsilon_s, id_H*S = \varepsilon_t, S*id_H*S = S,$$ where $*$ is the convolution product.

We remark that a weak Hopf algebra is a Hopf algebra if and only if $\Delta(1) = 1\otimes 1$ if and only if $\varepsilon$ is a homomorphism of algebras.

Let $H$ be a weak Hopf algebra and $A$ be an algebra. By [6], $A$ is a right weak $H$-comodule algebra if there is a right $H$-comodule structure $\rho_A: A \rightarrow A\otimes H, a\mapsto a_0\otimes a_1$ such that $\rho(ab) = \rho(a)\rho(b)$ for each $a, b \in A$ and $1_0\otimes1_1\in A\otimes H_t.$

The coinvariants defined by $$ B := A^{coH} := \{a\in A| \rho_A(a) = 1_0a\otimes 1_1\} $$ is a subalgebra of $A$. We say the extension $A/B$ is right weak $H$-Galois if the map $\beta: A\otimes_BA \rightarrow A\otimes_tH,$ given by $ a\otimes_Bb \mapsto ab_0\otimes_tb_1$ is bijective, where
$$A\otimes_tH:= (A\otimes H)\rho(1) = \{a1_0\otimes h1_1 \stackrel {\triangle}= a\otimes_t h|a\in A, h\in H\}.$$

 Recall from [3] that the finitistic dimension of $A$, is defined as

        fin.dim $A$ = sup$\{$pd$_A M\mid M\in A$-mod, pd$_A M < \infty\}.$

        The famous finitistic dimension conjecture says that fin.dim $A < \infty $ for any finite dimensional algebra $A$ and it has been proved
        that several classes of algebras have finite finitistic dimension. The finitistic dimension conjecture is still open now, for details, we refer to [11,12,14,16] and the references therein.

The purpose of this paper is to study the relationship of finitistic dimensions under right faithfully flat weak Hopf-Galois extension.
The main motivation for the study is that zhou and Li proved in [15]: let $H$ be a finite dimensional weak Hopf algebra and
$A/B$ be  a right faithfully flat weak $H$-Galois extension, if $H$ is semisimple, then the finitistic dimension of $A$ is less than or equal to that of $B$. Thus the question arises:   whether the finitistic dimension of $B$ is less or equal to that of $A$ when $A/B$ is  a right faithfully flat weak $H$-Galois extension? This is of particular interest because of the close relationship to the finitistic dimension conjecture.

         Our main result is as follows.

         \vskip 0.2in

{\bf Theorem.}  {\it Let $H$ be a finite-dimensional weak Hopf algebra over a field $k$ and $A/B$ be a right faithfully flat weak $H$-Galois extension. If the finitistic dimension of $B$ is finite, then it is less or equal to that of $A$.}

         \vskip 0.2in

In this paper, we always assume $H$ is finite dimensional and we follow the standard terminology and notation used in the representation theory of algebras and quantum groups, see [1,2,10].

\section { Proof of the Theorem}

\vskip 0.2in

\ \ \ \ \ \ The following well known Lemma shows the relations between the projective dimensions of the three modules in a short exact sequence, which will be used later.

\vskip 0.2in

{\bf Lemma 2.1.}  {\it If $0\rightarrow A_1 \rightarrow A_2 \rightarrow A_3 \rightarrow 0$ is exact and two of the modules have finite projective dimension, then so does the third. Moreover, if $n < \infty, {\rm pd}A_1 = n$, and {\rm pd}$A_3\leq n$, then {\rm pd}$A_2 = n$}.

Let $A/B$ be a right weak $H$-Galois extension. Consider the following two functors:

$A\otimes_B-: B-$Mod $\rightarrow A-$Mod, $M \mapsto A\otimes_B M,$

$_B(-): A-$Mod $\rightarrow B-$Mod, $M \mapsto  M,$

where $_B(-)$ is the restriction functor. It is known that $(A\otimes_B-, _B(-))$ is an adjoint pair. Let ($F, G$) be an adjoint pair of functors of abelian categories. If $F$ is exact, then $G$ preserves injective objects, if $G$ is exact, then $F$ preserves projective objects.

\vskip 0.2in

{\bf Lemma 2.2.}  {\it Let $A/B$ be a right weak $H$-Galois extension for a finite-dimensional weak Hopf algebra $H$. Then for each (finitely generated) B-module M,
{\rm pd}$_B(A\otimes_BM) \leq$ {\rm pd}$_A(A\otimes_BM ) \leq$ {\rm pd}$_BM$.}

\vskip 1mm

 {\bf Proof.} According to [4, Corallary 4.3], $A_B$ and $_BA$ are both finitely generated projective. It follows that pd$_B(A\otimes_BM) \leq$ pd$_A(A\otimes_BM) $ by the change of ring theorem.
 Assume that pd$_BM = n < \infty$, and let $\cal P$ be a projective resolution of $M$ as a $B$-module of length n. Now consider the adjoint pair $(A\otimes_B-, _B(-))$. Since $_B(-)$ is exact, the functor $A\otimes_B-$ preserves projective objects. It follows that  $A\otimes_B\cal P$ is a projective
 resolution of $A\otimes_BM$ as an $A$-module. This implies that pd$_A(A\otimes_BM) \leq $ pd$_BM$. The proof is completed.    \hfill$\Box$

 \vskip 0.2in

Now we can prove the main result of the paper.

 \vskip 0.2in

{\bf Theorem 2.3.}  {\it  Let $A/B$ be a right faithfully flat weak $H$-Galois extension for a finite-dimensional weak Hopf algebra $H$. If {\rm fin.dim} $B < \infty$, then {\rm fin.dim} $B\leq$ {\rm fin.dim} $A$.}

\vskip 1mm

 {\bf Proof.} For any left $B$-module $M$, there is a map $\varphi: M\rightarrow A\otimes_B M$ given by $m \mapsto 1\otimes m.$ Then it is a preliminary fact that $\varphi$ is monic since $A$ is a right faithfully flat $B$-module.

 Suppose fin.dim $B = n < \infty$, we choose $_BM$ with pd$_BM = n$. There is an $B$-exact sequence
 $$ 0 \rightarrow M \rightarrow A\otimes_B M \rightarrow C \rightarrow 0.$$
 By Lemma 2.2., pd$_B( A\otimes_B M) \leq$ pd$_BM $, so pd$_BC$ is also finite and  pd$_BC \leq n$ since fin.dim $B = n$. Now Lemma 2.1. gives
 pd$_B( A\otimes_B M) = n$.
By Lemma 2.2 again, since pd$_B(A\otimes_BM) \leq$ pd$_A(A\otimes_BM) \leq$ pd$_BM$, we conclude that pd$_A(A\otimes_BM) = n.$

Therefore, it follows directly that fin.dim $B\leq$ fin.dim $A$. The proof is completed.    \hfill$\Box$

 \vskip 0.2in

  {\bf Corollary 2.4.}  {\it  Let $A/B$ be a right faithfully flat weak $H$-Galois extension for a finite dimensional weak Hopf algebra $H$. Suppose $H$ is semisimple. If the finitistic dimension conjecture holds, then {\rm fin.dim} $B = $ {\rm fin.dim} $A$.}

  \vskip 1mm

  {\bf Proof.} According to [15], Let $A/B$ be a right faithfully flat weak $H$-Galois extension for a finite dimensional weak Hopf algebra $H$. Suppose $H$ is semisimple, then fin.dim $A \leq $ fin.dim $B$. By Theorem 2.3., the proof is completed.   \hfill$\Box$

 \vskip 0.2in

  Note that if $H$ is an ordinary Hopf algebra, then the weak $H$-Galois extension is just an $H$-Galois extension. So we may ask the following question:
Let $A/B$ be a right faithfully flat right $H$-Galois extension for a finite dimensional Hopf algebra $H$. Is the finitistic dimension of $B$  less or equal to that of $A$? Of course, An affirmative answer to the finitistic dimension conjecture would also give an affirmative answer to this question.

\begin{description}

\item{[1]}\ I.Assem, D.Simson, A.Skowronski, $Elements\ of\ the\ representation\ theory\ of\\ associative\ algebras$, Vol. 1, Cambridge
Univ. Press, 2006.

\item{[2]}\ M.Auslander, I.Reiten, S.O.Smal$\phi$, $Representation\ Theory\ of\ Artin\ Algebras$, Cambridge Univ. Press, 1995.

\item{[3]}\ H.Bass, Finitistic dimensions and a homological generalization of semiprimary rings, $Trans.\ Amer.\ Math.\ Soc$ 65(1960)466-488.

\item{[4]}\ G.B$\ddot{\rm o}$hm, $Galois theory for Hopf algebroids$, Ann. Univ. Ferrara Sez. VII(N.S.) 51(2005) 233-262.

\item{[5]}\ G.B$\ddot{\rm o}$hm, F.Nill, K.Szlach\'{a}nyi, Weak Hopf algebras I. Integral theory and C$^{*}$-structure, $J.Algebra$ 221(1999) 385-438.

\item{[6]}\ S.Caenepeel, E.De Groot, Modules over weak entwining structures, $Contemporary\\ Mathematics$ 267(2000) 31-54.

\item{[7]}\ T. Hayashi, Quantum group symmetry of partition functions of IRF models and its applications to Jones's index theory, $Comm.Math.Phys$ 157(1993) 331-345.

\item{[8]}\ L.Kadison, Galois theory for bialgebroids, depth two and normal Hopf subalgebras, $Ann. Univ.Ferrara Sez.VII(N.S.)$ 51(2005) 209-231.

\item{[9]}\ D.Nikshych, L.Vainerman, Finite quantum groupoids and their applications, in : New Directions in Hopf Algebras, $MSRI Publications$ 43(2002) 211-262.

\item{[10]}\ M.Sweedler, $Hopf Algebras$, W.A.Benjamin, New York, 1969.

\item{[11]}\ C.Wang, C.Xi,  Finitistic dimension conjecture and radical-power extensions, $J.Pure Appl.Algebra$221(4)(2017) 832-846.

\item{[12]}\ J.Wei, Finitistic dimension and Igusa-Todorov algebras, $Adv.Math$ 222(6)(2009) 2215-2226.

\item{[13]}\ T.Yamanouchi, Duality for generalized Kac algebras and a characterization of finite groupoid algebras, $J.Algebra$163(1994) 9-50.

\item{[14]}\ A.Zhang, S.Zhang, On the finitistic dimension conjecture of Artin algebras, $J.Algebra$ 320(2008) 252-258.

\item{[15]}\ X.Zhou, Q.Li, Finitistic dimensions of weak Hopf-Galois extensions, $Filomat$ 30,10(2016) 2825-2828.

\item{[16]}\ B.Zimmermann-Huisgen, The finitistic dimension conjectures-A tale of 3.5 decades, in: Abelian Groups and Modules, Padova,
1994, in: Math.Appl.,vol.343, Kluwer Academic Publishers, Dordrecht, (1995) 501-517.

\end{description}

\end{document}